\documentclass[11pt]{amsart}

\usepackage{a4wide,latexsym,amssymb,amsthm,amsmath}
\usepackage{graphicx,color}

\theoremstyle{plain}

\newtheorem{theorem}{Theorem}

\newtheorem{proposition}{Proposition}
\newtheorem{corollary}{Corollary}

\newtheorem{remark}{Remark}
\newtheorem{example}{Example}

\def\dx{\partial_x}
\def\dy{\partial_y}
\def\du{\partial_u}
\def\dv{\partial_v}

\def \n {\nabla}
\def \b {\beta}
\def \a {\alpha}

\def \grad {{\rm  grad}\ }

\def\<{ \left < }
\def\>{ \right > }

\def\R{\mathbb{R}}

\def\E{\mathbb{E}}
\def\H{\mathbb{H}}
\def\S{\mathbb{S}}

\def\<{\langle}
\def\>{\rangle}

\begin{document}

\title[]
{On certain surfaces in the Euclidean space ${\mathbb{E}}^3$}

\author[M.~I.~Munteanu]{Marian Ioan Munteanu}
\thanks{The first author was supported by Grant PN-II ID 398/2007-2010 (Romania)}
\email{marian.ioan.munteanu (at) gmail.com}
\author[A.~I.~Nistor]{Ana-Irina Nistor}
\email{ana.irina.nistor (at) gmail.com}

\begin{abstract}

In the present paper we classify all surfaces in $\E^3$ with a canonical principal direction.
Examples of these type of surfaces are constructed. We prove that the only minimal surface with a
canonical principal direction in the Euclidean space
${\mathbb{E}}^3$ is the catenoid.

\end{abstract}



\maketitle


\section{Preliminaries}

Due to recent research work in the field of classical differential geometry, the theory of surfaces and submanifolds knew a rapid development. Next to the classical problems of minimality
and flatness for different types of surfaces, another topic is represented by the study of constant angle surfaces. Even though in the Euclidean space the constant angle surfaces
are known in literature, in \cite{mn:MN09} is given a new approach of this problem regarding the ambient space $\E^3$ as the product space $\R^2\times\R$. By definition,
a constant angle surface is defined as a surface for which its unit normal makes a constant angle with a fixed direction given by the real line $\R$. Projecting the fixed direction
on the tangent plane to the surface and denoting by $U$ its tangent part we get that $U$ is a principal direction with null corresponding principal curvature. Assuming that $U$ remains a principal
direction but the corresponding principal curvature is different from zero -  the angle function is no longer constant - we denominate $U$ a canonical principal direction. First result on this topic
was given in \cite{mn:DFV09} for the ambient space $\S^2\times\R$ and the study in $\H^2\times\R$ was done in \cite{mn:DMN10}.

In the present paper we classify all surfaces with a canonical principal direction in $\E^3$. In our study we make use of canonical coordinates on the surface, obtaining also
classification theorems under the extra assumptions of minimality or flatness.
For example, we prove that the only minimal surface with a canonical principal direction in the Euclidean space
${\mathbb{E}}^3$ is the catenoid and we give its parametrization in canonical coordinates.
Moreover, illustrative examples of angle functions are constructed for known surfaces in $\E^3$ under harmonicity restrictions.

\smallskip

Let us consider a surface $M$ isometrically immersed in
$\E^3$ endowed with the scalar product $\<~,~\>$ and with the flat connection $\widetilde{\nabla}$. Denote by $g$ the
metric on $M$ which is the restriction of the scalar
product on $M$ and by $\nabla$ its corresponding Levi-Civita
connection. We consider an orientation of $\E^3$ and we denote by
$\overrightarrow{k}$ the fixed direction. If $N$ represents the unit
normal to the surface, then
$\theta(p):=\angle{(N,\overrightarrow{k})}$, with $\theta(p)\in(0,\pi)$,
represents the angle function between the unit normal and the fixed
direction in any point of the surface $p\in M$.

Classically we have the Gauss and Weingarten formulas for the surface $M$ isometrically immersed in $\E^3$:

{\bf(G)}\qquad\qquad $\widetilde{\nabla}_X Y=\nabla_X Y + h(X,Y)$

{\bf(W)}\qquad\qquad $\widetilde{\nabla}_X N=-A X$

for every $X,Y$ tangent to $M$. Moreover $h$ is a symmetric (1, 2)-tensor field
called the second fundamental form of the surface and $A$ is
a symmetric (1, 1)-tensor field denoting the shape operator
associated to $N$ which satisfies
$\<h(X,Y),\ N\>=g(X,AY)$ for any vector fields $X,Y$ tangent to $M$.

Denoting by $R$ the curvature tensor on $M$ and using the previous
notations, the equations of Gauss and Codazzi are given by

$
{\bf (E.G.)}\qquad\qquad
R(X,Y)=AX\wedge AY
$

$
{\bf (E.C.)} \qquad\qquad
(\nabla_X A)Y-(\nabla_Y A)X=0
$

where $X\wedge Y\in{\mathcal{T}}^1_1(M)$, $\big(X\wedge Y\big)Z=g(X,Z)Y-g(Y,Z)X$, for all $X,Y$ tangent to $M$.

One can decompose the fixed direction $\overrightarrow{k}$ as
\begin{equation}
\label{mn:k}
\overrightarrow{k}=U(p)+\cos\theta(p)N(p),
\end{equation}
where $U(p)$ is a tangent vector to $M$ in a point $p$ of the surface.

It follows that $\cos\theta(p)=\<\overrightarrow{k},N(p)\>$. From now on we drop the
explicit writing of the argument $p$ being obvious that the relations involving the angle function
are local and take place
in a neighborhood of any point of the surface.
All objects we use in this paper are supposed to be smooth, at least locally.
Moreover, $\theta\neq 0$ and $\theta\neq\frac{\pi}{2}$ because these situations were already studied
as particular cases of constant angle surfaces \cite{mn:MN09}.

Taking into account the decomposition \eqref{mn:k}, from the equation {\bf (E.G.)} the Gaussian curvature
can be computed as
\begin{equation}
\label{mn:K}
K=\det A.
\end{equation}

\begin{proposition}
\label{mn:p1}
For any $X$ tangent to $M$ the following statements hold
\begin{eqnarray}
\label{mn:p1_1}
& & \nabla_X U = \cos\theta AX \\
\label{mn:p1_2}
& & X[\cos\theta]=-g(AU,\ X),\ {\emph which\ can\ be\ equivalently\ written\ as} \\
\label{mn:p1_3}
& & AU=\sin\theta\ \grad\theta.
\end{eqnarray}
Here $\grad$ is thought with respect to the metric $g$.
\end{proposition}

\begin{proof}
Let us compute $\widetilde{\nabla}_X \overrightarrow{k}=\nabla_X U + h(X,\ U)+X[\cos\theta]N - \cos\theta AX$.
Identifying the tangent parts and taking into account that
$\widetilde{\nabla}_X \overrightarrow{k}=0$ we get \eqref{mn:p1_1}.
Next, identifying the normal parts and using the fact that $h(X,\ U)=g(AX,\ U)N$ one has \eqref{mn:p1_2}.
In order to retrieve \eqref{mn:p1_3} we can write $X[\cos\theta]=-\sin\theta~d\theta(X)$
and together with \eqref{mn:p1_2} yields $\displaystyle AU=\sin\theta~d\theta^{\sharp}$, where $\sharp$ denotes
the rising indices operation with respect to the metric $g$. \
We conclude with $g(d\theta^{\sharp},\ X)= X(\theta)=g({\rm grad\ }\theta,\ X)$.

\end{proof}

In the sequel we propose a way to deal with orthogonal coordinates on $M$.

\begin{proposition}
\label{mn:p2} For any angle function
$\theta\notin\{0,\frac{\pi}{2}\}$ one can choose local
coordinates $(x,y)$ on the surface $M$, isometrically immersed in $\E^3$, with
$\partial_x$ in direction of $U$ and such that the metric has the form
\begin{equation}
\label{mn:p2_1}
g=\frac{1}{\sin^2\theta}dx^2+\b^2(x,y)dy^2.
\end{equation}
The shape operator in the basis $\{\partial_x,\ \partial_y\}$ can be
expressed as
\begin{equation}
\label{mn:p2_2}
A = \left( \begin{array}{ccc}
\theta_x\sin\theta &  \theta_y\sin\theta\\ [2mm]
\frac{\theta_y}{\sin\theta\b^2} & \frac{\sin^2\theta\b_x}{\cos\theta\b}\end{array}
\right)
\end{equation}
and the functions $\b$ and $\theta$ are related by the PDE
\begin{equation}
\label{mn:p2_3} \frac{\sin^2\theta}{\cos\theta}\frac{\b_{xx}}{\b} +
\frac{\sin\theta\theta_x}{\cos^2\theta}\frac{\b_{x}}{\b}
+\frac{\theta_{y}}{\sin\theta}\frac{\b_y}{\b^3}+
\left(2\frac{\cos\theta\theta_y^2}{\sin^2\theta}-\frac{\theta_{yy}}{\sin\theta}\right)\frac{1}{\b^2}=0.
\end{equation}
\end{proposition}

\begin{proof}
Choosing an arbitrary point $p\in M$ such that the angle function
$\theta\neq 0,\frac{\pi}{2}$ we can consider locally the orthogonal
coordinates $(x,y)$ such that $\partial_x$ is in direction of $U$
and the metric is given by
\begin{equation}
\label{mn:pp2_11}
g=\alpha^2(x,\ y)dx^2+\beta^2(x,\ y)dy^2
\end{equation}
with $\a$ and $\b$ functions on $M$. The Levi-Civita connection for
this metric can be expressed in terms of $(x,y)$ coordinates as follows:
\begin{subequations}
\renewcommand{\theequation}{\theparentequation .\alph{equation}}
\label{mn:pp2}
\begin{eqnarray}
&& \label{mn:pp2_1}
\n_{\dx} \dx =\frac{\a_x}{\a}\dx - \frac{\a\a_y}{\b^2}\dy\\[2mm]
&& \label{mn:pp2_2}
\n_{\dx} \dy = \n_{\dy} \dx = \frac{\a_y}{\a}\dx + \frac{\b_x}{\b}\dy\\[2mm]
&&
\label{mn:pp2_3}
\n_{\dy} \dy =-\frac{\b\b_x}{\a^2}\dx + \frac{\b_y}{\b}\dy
\end{eqnarray}
\end{subequations}
One may compute the shape operator $A$ in this way. Since $\partial_x$ is in the direction of $U$, then
$\displaystyle U=\frac{\sin\theta}{\a}\partial_x$. Combining now the
expression $\displaystyle
\grad\theta=\frac{\theta_x}{\alpha^2}\dx+\frac{\theta_y}{\b^2}\partial_y$
with \eqref{mn:p1_3} we get
\begin{equation}
\label{mn:pp2_4}
A\partial_x=\frac{\theta_x}{\a}\partial_x+\frac{\theta_y\alpha}{\b^2}\partial_y.
\end{equation}
On the other hand, computing $A\partial_x$ using formulas \eqref{mn:p1_1} for $X=\dx$ and \eqref{mn:pp2_1}
we have
\begin{equation}
\label{mn:pp2_7}
A\dx = \frac{\theta_x}{\a}\dx - \tan\theta\frac{\a_y}{\b^2}\dy.
\end{equation}
Comparing \eqref{mn:pp2_4} and \eqref{mn:pp2_7} it follows that $\theta$ and $\a$ are related by
\begin{equation}
\label{mn:pp2_6}
\tan\theta\a_y+\a\theta_y=0.
\end{equation}

Moreover, in order to determine $A\partial_y$, we use formulas
\eqref{mn:p1_1} for $X=\dy$ and \eqref{mn:pp2_2} getting
\begin{equation}
\label{mn:pp2_5}
A\dy = \frac{\theta_y}{\a}\dx + \tan\theta\frac{\b_x}{\a\b}\dy.
\end{equation}

Hence, the shape operator is given by
\begin{equation}
\label{mn:pp2_8}
A = \left(
\begin{array}{ccc}
\frac{\theta_x}{\a} & \frac{\theta_y}{\a}  \vspace{2mm} \\
\frac{\theta_y\a}{\b^2} & \frac{\tan\theta\b_x}{\a\b}
\end{array}\right).
\end{equation}

The expression \eqref{mn:pp2_6} is equivalent with
$\dy(\a\sin\theta)=0$ and it yields $\displaystyle
\a=\frac{\phi(x)}{\sin\theta}$, where $\phi$ is a function on $M$
depending on $x$. Changing the $x$-coordinate we can assume that
$\displaystyle \a=\frac{1}{\sin\theta}$ and substituting it in the
general expression of the metric \eqref{mn:pp2_11} we get
\eqref{mn:p2_1}. Moreover, replacing the value of $\a$ in expression \eqref{mn:pp2_8}
of $A$ we obtain the shape operator given
exactly by formula \eqref{mn:p2_2}.

Furthermore, {\bf (E.C.)} is equivalent with $\nabla_{\dx}
A\dy-\nabla_{\dy}A\dx=0$. By straightforward computations the PDE
\eqref{mn:p2_3} is obtained, concluding the proof.

\end{proof}

\begin{remark}\rm
Any two functions $\theta$ and $\b$ defined on
a smooth simply connected surface $M$ related by \eqref{mn:p2_3} give
a surface isometrically immersed in $\E^3$ with the metric in the form \eqref{mn:p2_1} and the shape operator
 \eqref{mn:p2_2}.
\end{remark}

{\em Sketch of Proof.} Knowing $\theta$ and $\beta$ one could write the metric of the surface in form \eqref{mn:p2_1} and one
could determine the coefficients of the second fundamental form such that the associated matrix in the $\{\dx,\ \dy\}$
basis is given by \eqref{mn:p2_2}. The existence of the immersion easily follows
applying the Fundamental Theorem for the local theory of surfaces and Proposition~\ref{mn:p2}.

\hspace{152mm}\raisebox{-0.2ex}{\Large$\Box$}
\smallskip

As we would like to find some explicit parameterizations,
we should be able to solve \eqref{mn:p2_3} in order determine the
metric. A first step to solve it is to impose some extra conditions getting some
results involving harmonic maps and illustrative examples.

\begin{proposition}
\label{mn:p3}
Let $M$ be a minimal isometric immersion in $\E^3$.
We can choose local coordinates $(x,y)$ on $M$  such that $\partial_x$ is in direction of $U$,
the metric of the surface can be expressed as
\begin{equation}
\label{mn:p3_1}
g=\frac{1}{\sin^2\theta}(dx^2+dy^2)
\end{equation}
and the shape operator $A$ in the basis $\{\dx,\ \dy\}$ has the following expression
\begin{equation}
\label{mn:p3_2}
A = \sin\theta\left( \begin{array}{ccc}
\theta_x & \theta_y\\
\theta_y&-\theta_x
\end{array} \right).
\end{equation}
Moreover, the function $\displaystyle \log\left(\tan\frac{\theta}{2}\right)$ is harmonic.
\end{proposition}
\begin{proof}
Using the results from \emph{Proposition~$\ref{mn:p2}$}, the minimality condition
${\rm trace} A=0$ in \eqref{mn:p2_2}
yields the PDE
$\cos\theta\theta_x\b+\sin\theta\b_x=0$, or equivalently $(\b\sin\theta)_x=0$.
Integrating once w.r.t $x$ and after a change of $y$-coordinate, one finds
$\b=\frac{1}{\sin\theta}$.
Hence, combining it with \eqref{mn:p2_1} we get the metric \eqref{mn:p3_1},
which corresponds to isothermal coordinates $(x,y)$ on $M$.
The expression of the shape operator \eqref{mn:p3_2} follows after straightforward computations.

Now, condition \eqref{mn:p2_3} together with the expression of $\b$ yields
\begin{equation}
\label{mn:pp3_1}
\cos\theta(\theta_x^2+\theta_y^2)-\sin\theta(\theta_{xx}+\theta_{yy})=0.
\end{equation}
The Laplacian of the surface $M$ is $\Delta=\sin^2\theta(\partial_{xx}^2+\partial_{yy}^2)$,
and therefore the above equation is equivalent to
\begin{equation}
\label{mn:pp3_2}
\Delta\log\left(\tan\frac{\theta}{2}\right)=0.
\end{equation}
Hence, the considered function is harmonic and this concludes the proof.

\end{proof}
\begin{corollary}
There are no minimal, compact and orientable surfaces isometrically immersed in $\E^3$.
\end{corollary}
{\em Sketch of proof.} We proceed by contradiction. If $M$ is such a surface (minimal, compact and orientable)
and we denote by $\theta$ the angle function then we could apply the previous proposition obtaining that
$\log \left(\tan\frac{\theta}{2}\right)$ is harmonic. By compactness of $M$ it follows that $\theta$ is constant (see e.g. \cite{mn:Mo01}).
Accordingly to the classification given  in \cite{mn:MN09} we get the contradiction.
(All constant angle surfaces in $\E^3$ are ruled surfaces, hence they cannot be compact.)

\hspace{152mm}\raisebox{-0.2ex}{\Large$\Box$}

\begin{remark}\rm
Any smooth function $\theta$ defined on a smooth simply connected surface $M$ satisfying
\eqref{mn:pp3_2} gives a minimal surface $M$ in $\E^3$ such that the metric on the surface
can be written in the form \eqref{mn:p2_1} and the shape operator is given by \eqref{mn:p2_2}.
\end{remark}

\smallskip

At this point we are interested to give some examples of angle functions $\theta$ for which
the corresponding surface is minimal in $\E^3$.
So, we have to solve \eqref{mn:pp3_1}.

In order to do this, let us look for $\theta$ such that there exists a real constant $a$ satisfying
$\theta_x=a\theta_y$. Computing $\theta_{xx}=a\theta_{xy}$ and
$\theta_{yy}=\frac{\theta_{xy}}{a}$, equation \eqref{mn:pp3_1} becomes
$
\cos\theta\theta_y^2-\sin\theta\frac{\theta_{xy}}{a}=0.
$
This yields $\partial_x\Big(\frac{\theta_y}{\sin\theta}\Big)=0$ and
$\partial_y\Big(\frac{\theta_y}{\sin\theta}\Big)=0$, which means
$\displaystyle\frac{\theta_y}{\sin\theta}=b,\ b\in\R$. Integrating now with respect to $y$ it follows
$\displaystyle\ln\left|\tan\frac{\theta}{2}\right|=by+m(x)$, where $m$ is a function depending on $x$
which must be determined. Taking the derivative with respect to $x$ in the previous expression,
we find $m(x)=abx+d$, with $d\in \R$. After simple computations one concludes that
\begin{equation}
\label{mn:pp3_7}
\theta=2\arctan\left(c~e^{b(ax+y)}\right),\
{\rm where}\ b,c\in \R\ .
\end{equation}

Hence, there exists a minimal surface $M$ isometrically immersed in $\E^3$
for which the normal to surface forms with the fixed direction
$\overrightarrow{k}$ the angle $\theta$ given by \eqref{mn:pp3_7}.

\begin{remark}\rm
Notice that $b(ax+y)+d$ is a harmonic function.
We ask now if for any harmonic function generically denoted $f$ the angle function
\begin{equation}
\label{mn:pp3_9}
\theta=2\arctan\left(e^{f}\right)
\end{equation}
still gives a minimal surface ?
\end{remark}

\begin{proof}
Indeed, the answer is positive, from Proposition~\ref{mn:p3}. Yet,
in order to see this, suffices to check that $\theta$ fulfills \eqref{mn:pp3_2},
which is equivalent in the statement of Proposition~\ref{mn:p3}
with the fact that function $\displaystyle\log\left(\tan\frac{\theta}{2}\right)$ is harmonic.
So, for the expression of $\theta$ given by \eqref{mn:pp3_9} we compute

$\displaystyle \sin\theta=2\frac{e^f}{1+e^{2f}}
\qquad\qquad\qquad\qquad\qquad\qquad\qquad \cos\theta=\frac{1-e^{2f}}{1+e^{2f}}$

$\displaystyle \theta_x=2\frac{e^f f_x}{1+e^{2f}}
\qquad\qquad\qquad\quad\qquad\qquad\qquad\qquad \theta_y=2\frac{e^f f_y}{1+e^{2f}}$

$\displaystyle \theta_{xx}=2\frac{e^f(1+e^{2f})(f_x^2+f_{xx})-2e^{3f}f_x^2}{(1+e^{2f})^2}\quad\qquad\quad
\theta_{yy}=2\frac{e^f(1+e^{2f})(f_y^2+f_{yy})-2e^{3f}f_y^2}{(1+e^{2f})^2}.
$

Using these expressions and taking into account that $f$ is harmonic,
namely $f_{xx}+f_{yy}=0$, we get that \eqref{mn:pp3_1} is automatically
satisfied. Hence, in order to give more examples of minimal surfaces
we use Proposition~\ref{mn:p3} with the angle function given by \eqref{mn:pp3_9}
for any $f$ - harmonic function on an open set of $\R^2$.

\end{proof}

\begin{example}\rm
A first example of angle $\theta$ that corresponds to a minimal surface in $\E^3$,
follow-on the previous remark, can be obtained taking the harmonic function
$f:\R^2\setminus\{0\}\longrightarrow\R^2,\  f(x,\ y)=\ln(x^2+y^2)$
which, in physics, expresses the electric potential due to a line charge,
for which the angle function that determines the surface is given by
$\theta=2\arctan(x^2+y^2)$.
\end{example}

Thinking now conversely, we believe that it is interesting to study the angle function $\theta$
and to find corresponding isothermal coordinates
in order that Proposition~\ref{mn:p3} holds true for the well-known examples of minimal surfaces in $\R^3$,
namely the {\em helicoid}, the {\em catenoid}, the {\em Enneper surface} and the
{\em Scherk surface} respectively.
We will start for each of them with the usual parametrization.

\begin{example}\rm \label{mn:helicoid}
The classical parametrization for helicoid, denoted by $\mathcal{H}$, is given by
$$
r(u,v)=(u\cos v,\ u\sin v,\ v).
$$
The objects we are looking for are: $g=du^2+(u^2+1)dv^2$ and $\theta=2\arctan(\sqrt{u^2+1}-u)$.

It can be checked that $\Delta\log\left(\tan\frac{\theta}{2}\right)=0$,
where $\Delta=\du^2+\frac{1}{u^2+1}\dv^2+\frac{u}{u^2+1}\du$ is the Laplacian on the surface $\mathcal H$.
Since the above coordinates $(u,v)$ are not isothermal,
we make the change of coordinates $u=\sinh x$ and $v=y$ such that
the metric becomes $g=\cosh^2 x (dx^2+dy^2)$ and we stay in the hypothesis of Proposition~\ref{mn:p3}.
\end{example}

\begin{example}\rm
Studying the catenoid $\mathcal C$ parameterized by
$$
r(u,v)=(\cosh u\cos v,\ \cosh u\sin v,\ u)
$$
we get the metric $g=\cosh^2 u(du^2+dv^2)$ and we find the angle function
$\theta=2\arctan e^{-u}$. Since the Laplacian on $\mathcal C$ is
$\Delta=\frac{1}{\cosh^2v}(\partial_{uu}^2+\partial_{vv}^2)$
one easily obtains the harmonicity.
We also remark that $(u,v)$ are isothermal coordinates, hence Proposition~\ref{mn:p3} is again verified.
\end{example}

\begin{example}\rm
The parametrization of the Enneper surface $\mathcal E$ is given by
$$
r(u,v)=\Big(u-\frac1 3\ u^3+uv^2,\ -v+\frac1 3\ v^3-u^2v,\ u^2-v^2\Big),\ (u,v)\neq (0,0).
$$
The metric has the form $g=(1+u^2+v^2)^2(du^2+dv^2)$
and the angle function is\\ $\theta=2\arctan\frac{1}{\sqrt{u^2+v^2}}\ $.
Again, using the expression of the Laplacian
$\Delta=\frac{1}{(1+u^2+v^2)^2}(\partial{uu}^2+\partial_{vv}^2)$ on the surface
$\mathcal E$ one obtains that the function $\log\left(\tan\frac\theta 2\right)$ is harmonic. Moreover,
the coordinates $(u,v)$ are isothermal as in Proposition~\ref{mn:p3}.
\end{example}

\begin{example}\rm
The parametrization of the Scherk surface $\mathcal S$ over the square (see \cite{mn:ON06}) can be written as
$$
r(u,v)=\Big(u,\ v,\ \log\frac{\cos u}{\cos v}\Big).
$$
The angle function $\theta$ satisfies
$ \cos\theta = \Big(\frac{1}{\cos^2u}+\frac{1}{\cos^2v}-1\Big)^{-\frac1 2}$. Notice that
in this case the coordinates are no longer orthogonal, since the metric has the following form
$g=\frac{1}{\cos^2 u}du^2-2\frac{\sin u\sin v}{\cos u\cos v}dudv+\frac{1}{\cos^2v}dv^2$.
Looking for an isothermal parametrization in $(x,y)$ coordinates, one has to find $u=u(x,y)$ and $v=v(x,y)$
such that the following system is fulfilled:
$$
\left\{
\begin{array}{l}
 (u_x^2-u_y^2)\frac{1}{\cos^2u}+2(-u_x v_x +u_y v_y)\tan u\tan v +(v_x^2-v_y^2)\frac{1}{\cos^2 v}=0\\ [3mm]
 u_x u_y\frac{1}{\cos^2u}-(u_x v_y+u_yv_x)\tan u \tan v +v_xv_y\frac{1}{\cos^2 v}=0.
\end{array}
\right.
$$
The isothermal parametrization of $S$ is the given by
$r(x,y)=\Big(u(x,y),\ v(x,y),\ \log{\frac{\cos u(x,y)}{\cos v(x,y)}}\Big)$ with
$u(x,y)=\arctan\frac{2x}{1-x^2-y^2}$ and
$v(x,y)=\arctan\frac{-2y}{1-x^2-y^2}$ .
\end{example}

We conclude this section with a non existence result
\begin{remark}
There are no minimal and flat surfaces isometrically immersed in $\E^3$ with a non constant angle function.
\end{remark}
\begin{proof}
Computing the Gaussian curvature $K$ for a minimal surface given as in
Proposition~\ref{mn:p3} we get that $K=0$ is equivalent with
$\theta_x^2+\theta_y^2=0$. Consequently, it follows that $\theta$ is constant.

\end{proof}

\section{Surfaces with a Canonical Principal Direction}

The study of constant angle surfaces in $\E^3$ can be generalized for surfaces
whose angle function is no longer constant,
but certain properties are preserved. More precisely, this is the case in which $U$ remains a principal direction,
whereas the corresponding principal curvature is different from $0$. They will be called
\emph{surfaces with a canonical principal
direction}. We characterize these surfaces in the following

\begin{theorem}
Let $M$ be an isometrically immersed surface in $\E^3$. Let $(x,y)$ be orthogonal coordinates on $M$ such that $U$ is
collinear to $\dx$. Then, $U$ is a principal direction on $M$ everywhere if and only if $\theta_y=0$.
\end{theorem}

\begin{proof}
We know that such coordinates exist as in the proof of {\it Proposition}~\ref{mn:p2}.
We have
$$
U=\sin^2\theta\dx.
$$
Moreover, from the expression \eqref{mn:p2_2} of the shape operator it follows that
$$
AU=\sin^3\theta\theta_x\dx+\sin\theta\frac{\theta_y}{\b^2}\dy.
$$
We find that $U$ is a principal direction implies $\theta_y=0$.

Conversely, from \eqref{mn:p1_2} it follows $g(AU,\dy)=0$ which
means that $AU$ is parallel to $\dx$, hence $U$ is a principal direction for $M$.

\end{proof}

The following statement is essential for the rest of the paper.
\begin{proposition}
\label{mn:p4}
Let $M$ be a surface immersed in $\E^3$ and a point $p\in M$ such that
$\theta(p)\notin \{0,\frac{\pi}{2}\}$. If $U$ is a principal direction of $M$,
we can choose coordinates $(x,\ y)$ in a neighborhood of $p$ such that $\dx$ is in the direction of $U$,
the metric has the form
\begin{equation}
\label{mn:p4_1}
g=dx^2+\b^2(x,y)dy^2
\end{equation}
and the shape operator is given by
\begin{equation}
\label{mn:p4_2}
A = \left( \begin{array}{ccc}
\theta_x & 0 \\
0&\tan\theta\frac{\b_x}{\b}
\end{array} \right).
\end{equation}
Moreover, $\theta$ and $\b$ are related by the PDE
\begin{equation}
\label{mn:p4_3}
\b_{xx} + \tan\theta\theta_x\b_{x}=0
\end{equation}
and $\theta_y=0$.
\end{proposition}
\begin{proof}
 The results are obtained using similar techniques as in Proposition~\ref{mn:p2} by straightforward computations.

\end{proof}

\begin{remark}\rm
Accordingly to \cite{mn:DMN10}, we say that $(x,y)$ are {\em canonical coordinates} on $M$
if $U$ is a principal direction collinear to $\dx$ and the metric $g$ has the form
\eqref{mn:p4_1}. Notice that such coordinates are not unique and
two pairs $(x,y)$ and $(\overline{x},\overline{y})$ of canonical coordinates
are related by $\overline{x}=\pm\ x+c,\ c\in\R$ and $\overline{y}=\overline{y}(y)$.
\end{remark}

An illustration of canonical coordinates is given by {\bf Example~\ref{mn:helicoid}}, the classical parametrization
of the helicoid. In this case the metric is written in form \eqref{mn:p4_1} with $\b=\sqrt{u^2+1}$ and together
with $\theta=2\arctan(\sqrt{u^2+1}-u)$ fulfill \eqref{mn:p4_3} identically.

In order to determine explicitly all surfaces in $\E^3$ with a canonical principal
direction we have to solve \eqref{mn:p4_3}
in order to find the unknown function $\b$ from the expression \eqref{mn:p4_1} of the metric.
They are described in the following classification theorem:
\begin{theorem}\label{mn:t2}
A surface $M$ isometrically immersed in $\E^3$ with $U$ a canonical principal direction is given
(up to isometries of $\E^3$) by one of
the following cases:
\begin{itemize}
\item {\bf Case 1.}
\begin{equation}
\label{mn:t2_1}
r:M\rightarrow\E^3,\ r(x,\ y)=\left(\phi(x)(\cos y,\ \sin y)+\gamma(y),\ \int_0^x\sin\theta(\tau)d\tau\right)
\end{equation}
where
$$
\gamma(y)=\left(-\int_0^y \psi(\tau)\sin\tau d\tau,\ \int_0^y \psi(\tau)\cos\tau d\tau\right)
$$
\item{\bf Case 2.}
\begin{equation}
\label{mn:t2_2}
r:M\rightarrow\E^3,\ r(x,\ y)=\left(\phi(x)\cos y_0,\ \phi(x)\sin y_0,\
\int_0^x\sin\theta(\tau)d\tau\right)+y{\mathrm{v_0}}
\end{equation}
where
${\mathrm{v_0}}=(-\sin y_0,\ \cos y_0,\ 0),\ y_0\in\R$.
Notice that these surfaces are cylinders.
\end{itemize}
In both cases $\phi(x)$ denotes a primitive of $\cos\theta$.
\end{theorem}
\begin{proof}
Let us denote the isometric immersion of the surface $M$ in $\E^3$ by
$$
r:M\rightarrow\E^3,\ r(x,y)=\Big(r_1(x,y),\ r_2(x,y),\ r_3(x,y)\Big)=\Big(r_j(x,y),\ r_3(x,y)\Big),\ j=1,2.
$$
Since the statements of Proposition~\ref{mn:p4} hold true, we are able to choose canonical coordinates $(x,y)$
such that the
metric is given by \eqref{mn:p4_1}. At this point we have to determine the function $\b$ which satisfies the PDE
\eqref{mn:p4_3}, or equivalently,
$\partial_x\left(\frac{\b_x}{\cos\theta}\right)=0$. Integrating twice one gets:
\begin{itemize}
\item either $\b=k(y)(\phi(x)+\psi(y))$, where $\phi'(x)=\cos\theta$, $\psi(y)$ and $k(y)$ are defined on $M$

\item or $\b=\b(y)$. 
\end{itemize}

We may immediately obtain the $3^{\rm rd}$ component of the immersion $r$. Since from Proposition~\ref{mn:p4}
$U=\sin\theta\dx$, the decomposition \eqref{mn:k} becomes
\begin{equation}
\label{mn:tt2_8}
\overrightarrow k=\sin\theta r_x +\cos\theta N.
\end{equation}

Computing now
$(r_3)_x=\<r_x,\overrightarrow{k}\>=\sin\theta$ and
$(r_3)_y=\<r_y,\overrightarrow{k}\>=0$
we conclude that
\begin{equation}
\label{mn:tt2_9}
\displaystyle r_3(x,y)=\int^x\sin\theta(\tau)d\tau.
\end{equation}

From \eqref{mn:tt2_8} and \eqref{mn:tt2_9} one can express the normal to the surface as
\begin{equation}
\label{mn:tt2_10}
N=\Big(-\tan\theta(r_j)_x, \cos\theta\Big),\ j=1,2.
\end{equation}

Let us distinguish two cases for $\b$.

\textbf{Case 1.} $\b=k(y)(\phi(x)+\psi(y))$\\
After a change of the $y-$coordinate we may assume that
$
\b=\phi(x)+\psi(y)
$
and substituting it in \eqref{mn:p3_1} we get the metric
\begin{equation}
\label{mn:tt2_1}
g=dx^2+(\phi(x)+\psi(y))^2dy^2.
\end{equation}
By using Koszul formula one obtains the corresponding Levi-Civita connection
\begin{eqnarray*}
\label{mn:tt2_2} &&\n_{\dx} \dx =0,\ \ \ \ \n_{\dx} \dy = \n_{\dy} \dx =\frac{\cos\theta}{\phi(x)+\psi(y)}\dy \\
&&\n_{\dy} \dy =-(\phi(x)+\psi(y))\cos\theta\dx + \frac{\psi'(y)}{\phi(x)+\psi(y)}\dy.
\end{eqnarray*}

Taking into account the expression of the shape operator \eqref{mn:p4_2} and the metric \eqref{mn:tt2_1} we get
\begin{equation*}
\label{mn:tt2_5}
A = \left( \begin{array}{ccc}
\theta_x & 0 \\
0&\frac{\sin\theta}{\phi(x)+\psi(y)}
\end{array} \right).
\end{equation*}
Moreover, from the Weingarten formula {\bf (W)} we have
\begin{subequations}
\renewcommand{\theequation}{\theparentequation .\alph{equation}}
\label{eq:tt2_W}
\begin{eqnarray}
\label{mn:tt2_6}
&&N_x=-\theta_x\ r_x\\
\label{mn:tt2_7}
&&N_y=-\frac{\sin\theta}{\phi(x)+\psi(y)}\ r_y.
\end{eqnarray}
\end{subequations}
Computing the derivative with respect to $x$ in \eqref{mn:tt2_10},
$$
N_x=\left(-\frac{\theta_x}{\cos^2\theta}{(r_j)}_x -\tan\theta{(r_j)}_{xx}, -\sin\theta\theta_x\right)
$$
and combining it with \eqref{mn:tt2_6} it follows that $r$ must fulfil
$
\cos\theta{(r_j)}_{xx}+\sin\theta\theta_x {(r_j)}_x=0
$
which can be equivalently written $\displaystyle \partial_x \left(\frac{{(r_j)}_x}{\cos\theta}\right)=0$,
$j=1,2$. One gets
\begin{equation}
\label{mn:tt2_11}
\big(r_1,r_2\big)_x=\cos\theta f(y)
\end{equation}
where $f(y)=(\cos\varphi(y),\sin\varphi(y))$ represents a parametrization of the unit circle $\S^1$.
This is a consequence of the fact that $\|r_x\|^2=1$ which, combined with \eqref{mn:tt2_9} and \eqref{mn:tt2_11},
leads to $\|f(y)\|=1$. Here and all over this paper $\|\cdot \|$ denotes the Euclidean norm.

Integrating with respect to $x$ in \eqref{mn:tt2_11}  and taking into account \eqref{mn:tt2_9} we get
the following expression for the immersion $r$
\begin{equation}
\label{mn:tt2_12}
r(x,y) =\left(\phi(x) f(y)+\gamma(y),\ \int_0^x\sin\theta(\tau)d\tau\right)
\end{equation}
where $\gamma(y)=(\gamma_1(y), \gamma_2(y))$ is a smooth $\R^2-$valued map.

Since $r$ is an isometric immersion we get

\begin{itemize}
\item[(i)] $\<\gamma'(y),\ f(y)\>=0$\\
\item[(ii)]  $ \phi(x)^2\|f'(y)\|^2+\|\gamma'(y)\|^2+2\phi(x)\<f'(y),\ \gamma'(y)\>=\b^2$.
\end{itemize}

From (i) we deduce that $\gamma'(y)$ and $f'(y)$ are parallel vectors, so, there exists a $C^\infty$-function
$\eta(y)$  such that
$\displaystyle \gamma'(y)=\eta(y) f'(y)$. Replacing it in (ii)
and taking into account the expression of $\b$ we get equivalently
$\big(\phi(x)+\eta(y)\big)\left|\varphi'(y)\right|=\phi(x)+\psi(y)$.
Since $\phi(x)$ is not a constant, it follows $|\varphi'(y)|=1$ and consequently
$\eta(y)=\psi(y)$. By fixing an orientation on the $y$-axis and after a translation along it,
we may choose $\varphi(y)=y$ and hence
\begin{equation}
\label{mn:tt2_13}
\displaystyle \gamma(y)=\left(-\int_0^y \psi(\tau)\sin\tau d\tau,\ \int_0^y \psi(\tau)\cos\tau d\tau\right).
\end{equation}
Combining now \eqref{mn:tt2_9}, \eqref{mn:tt2_12} and \eqref{mn:tt2_13} we get exactly the parametrization
\eqref{mn:t2_1}.


\medskip

{\bf Case 2.} $\b=\b(y)$\\
After a change of the $y$-coordinate in this case, $\b=1$ and the metric becomes
$
g=dx^2+dy^2.
$
The shape operator is given by
$
A = \left( \begin{array}{cc}
\theta_x & 0 \\
0&0
\end{array} \right).
$
From the Weingarten formula {\bf (W)} we get
\begin{subequations}
\renewcommand{\theequation}{\theparentequation .\alph{equation}}
\label{eq:tt2_14_15}
\begin{eqnarray}
\label{mn:tt2_14}
&&N_x= - \theta_xr_x \\
\label{mn:tt2_15}
&&N_y=0.
\end{eqnarray}
\end{subequations}
Firstly, taking the derivative with respect to $y$ in \eqref{mn:tt2_10}
and combining it with \eqref{mn:tt2_15} one gets
\begin{equation}
\label{mn:tt2_16}
(r_j)_{xy}=0, \quad j=1,2.
\end{equation}
Secondly, taking the derivative with respect to $x$ in \eqref{mn:tt2_10} and combining it with \eqref{mn:tt2_14}
one obtains
\begin{equation*}
(r_j)_{xx}+\tan\theta\theta_x (r_j)_x=0, \quad j=1,2.
\end{equation*}
It follows
\begin{equation*}
\label{mn:tt2_17}
(r_j)_{x}=\cos\theta f_j(y),\ {\rm where}\  f_j\in C^{\infty}(M),\ j=1,2\ {\rm with}\ f_1(y)^2+f_2(y)^2=1.
\end{equation*}
Combining the above equation with \eqref{mn:tt2_16} one obtains
$$
(r_1(x,y),\ r_2(x,y))=f_{0}\phi(x)+\gamma(y)
$$ where
$f_{0}=(\cos y_0,\ \sin y_0)$ and $\gamma(y)=(-(\sin y_0)y+c_1, (\cos y_0)y +c_2),\ {\rm with }\ c_1,c_2\in \R$
since $r$ is an isometric immersion.
After a translation in $\E^3$, the constants $c_1$ and $c_2$ may be assumed to be zero.
So, the expression \eqref{mn:t2_2} is proven.

Conversely, we will give the proof only in {\bf Case 1} because the idea of the proof is the same also
in the second case.
Suppose that we have a surface given by \eqref{mn:t2_1} and we prove that it has $U$ as a
canonical principal direction. The tangent plane of the surface is generated by the following vectors
\begin{equation*}
\begin{array}{ll}
\hspace{-34mm} r_x=(\cos\theta\cos y,\ \cos\theta\sin y,\ \sin\theta)\vspace{2mm} \\
\hspace{-34mm} r_y=\left(-(\phi(x)+\psi(y))\sin y,\ (\phi(x)+\psi(y))\cos y,\ 0\right).
\end{array}
\end{equation*}
Hence, the metric is given by $g=dx^2+(\phi(x)+\psi(y))^2dy^2$, which corresponds to \eqref{mn:p4_1}.
By straightforward computations we get that the shape operator has the form as in \eqref{mn:p4_2} and from
its symmetry we have $\theta_y=0$.
One easily proves that $\<r_x,U\>=\sin\theta$ and $\<r_y,U\>=0$ concluding that $U$ is a principal
direction of the surface $M$ parameterized
by \eqref{mn:t2_1}. At this point the theorem is completely proved.

\end{proof}

An alternative proof of this result, but in a different manner, can be found in \cite{mn:To09}. The author classifies
hypersurfaces $f:M^n\rightarrow \mathbb{Q}_\epsilon^n\times\R$
with a principal direction, where $\mathbb{Q}_\epsilon^n$ denotes either the $n$-sphere $\S^n$,
the Euclidean $n$-space $\R^n$ or the hyperbolic $n$-space $\H^n$ according to
$\epsilon=1,\
\epsilon=0$, or $\epsilon=-1$.

The study of minimal surfaces is another classical problem in differential geometry.
Below we classify all minimal surfaces with a canonical
principal direction given by $U$ in $\E^3$.

\begin{theorem} 
Let $M$ be a surface isometrically immersed in $\E^3$. M is a minimal surface with $U$ a principal
direction if and only if the immersion is, up to
isometries of the ambient space, given by
\begin{equation}
\label{mn:t3_1}
r:M\rightarrow \E^3,\ \
r(x,\ y)=\left(\sqrt{x^2+c^2}\cos y,\ \sqrt{x^2+c^2}\sin y ,\ c\ln \big(x+\sqrt{x^2+c^2}\big)\right),\ c\in\R. \\
\end{equation}
\end{theorem}
\begin{proof}
The result is local, hence Proposition~\ref{mn:p4} can be used.
From the proof of the previous theorem we already know the solutions of \eqref{mn:p4_3},
namely $\b=\phi(x)+\psi(y)$, respectively $\b=1$ after a change of the $y$-coordinate.

Let us consider the first solution for $\b$.
Under the minimality assumption we get, after a translation along the $x$-coordinate, that
\begin{equation}
\label{mn:tt3_4}
\theta=\arctan\frac{c}{x}\ ,\ c\in \R.
\end{equation}
Writing $\beta$ in two ways, once $\beta=\frac{c}{\sin\theta}$ from the minimality condition
and second in the general form $\beta=\phi(x)+\psi(y)$, we immediately find that $\psi(y)$ is constant and it may
be added to $\phi(x)$. Therefore $\psi(y)$ could be considered zero.

Going back to Theorem~\ref{mn:t2} we know the parametrization \eqref{mn:t2_1}.
Substituting the value of $\theta$ from \eqref{mn:tt3_4}, by straightforward computations we find
$$
\phi(x)=\sqrt{x^2+c^2}\ , \quad
\gamma(y)=0
$$
and the $3^{\rm rd}$ component of the parametrization becomes
$$
\int^x\sin\theta(\tau)d\tau=c\ln(x+\sqrt{x^2+c^2}).
$$
Combining these expressions in \eqref{mn:t2_1}, we find that a minimal surface which has $U$ a
principal direction everywhere
is parameterized by \eqref{mn:t3_1}.

In the second case of the classification theorem, corresponding to $\b=1$, under the assumption of minimality
we get that $\theta_x=0$ which contradicts our initial hypothesis that $\theta$ is never constant.
So, this situation cannot occur.

The converse results immediately by direct computations.

\end{proof}

\begin{remark}\rm
We notice that this surface can be obtained by rotating the catenary around the $Oz$-axis.
Hence, the only minimal surface in the Euclidean 3-space with a canonical principal direction is
the catenoid. See for details \cite{mn:Bl75,mn:CD83}. This result can be retrieved in a different manner also from \cite{mn:To09}.
\end{remark}
\begin{theorem}
Let $M$ be a surface isometrically immersed in $\E^3$. M is a flat surface with $U$ a principal
direction if and only if the immersion is, up to
isometries of the ambient space, given by
\begin{equation*}
\label{mn:t4_2}
r:M\rightarrow \E^3,\quad
r(x,\ y)=\left(\phi(x)\cos y_0,\ \phi(x)\sin y_0,\ \int_0^x\sin\theta(\tau)d\tau\right)+ y{\rm v}_0
\end{equation*}
where
$
{\rm v}_0=\big(-\sin y_0, \cos y_0,0\big),\ y_0\in\R
$
and $\phi(x)$ represents a primitive of $\cos\theta$.
\end{theorem}
\begin{proof}
Using the canonical coordinates furnished by Proposition~\ref{mn:p4},
under the flatness assumption, from \eqref{mn:K} one gets $\theta_x\tan\theta\b_x=0$. As the angle function $\theta$
cannot be constant, it
yields $\b_x=0$ which implies $\b=\b(y)$. This corresponds precisely to {\bf Case 2} from the classification theorem.
Hence, the second case of the classification coincides with the class of all flat surfaces $M$ with a
canonical principal direction $U$.

\end{proof}

\begin{remark}\rm
Starting from Proposition~\ref{mn:p4} and using the expression of the shape operator $A$ together with
the expression \eqref{mn:p4_3} relating $\theta$ and $\beta$
we could find the expressions of the angle function $\theta$
in the case of equal principal curvatures or constant mean curvature (CMC).\linebreak
Studying surfaces with equal principal curvatures, we
solve \eqref{mn:p4_3} and after a change of $y$-co\-or\-di\-nate, we get the following solutions for $\beta$:

{\bf case 1.} $\b=\phi(x)+\psi(y)$ and replacing this expression in \eqref{mn:p4_2},
if the principal curvatures are equal, then $\theta$ must satisfy $\theta_x = \frac{\sin\theta}{\phi(x)+\psi(y)}$.
By direct computations and by taking into account that the angle function depends only of $x$ we get that
\begin{equation*}
\label{mn:eqcurv}
\theta=a x +b,\ a,b\in \R.
\end{equation*}
Note that in the case when the principal curvatures are constant, they must be equal with the same constant
and hence in this case we obtain a piece of the $2$-sphere
$\S^2$ in Euclidean $3$-space. Conversely, writing the sphere in canonical coordinates
(which coincide with the spherical coordinates) the angle function $\theta$ is
an affine function.

{\bf case 2.} $\beta=1$ implies $\theta_x=0$, hence $\theta$ is constant, situation excluded all over this paper.
\vspace{5mm}

Regarding the constant mean curvature surfaces, we discuss again the two solutions for $\beta$. \linebreak
{\bf case 1.} $\b=\phi(x)+\psi(y)$, then under the CMC$-$condition we get that $\theta$ must satisfy the following
differential equation:
$$
\theta_x + \frac{\sin\theta}{\phi(x)+\psi_0}=2H,\
{\rm where}\ H\ {\rm denotes\ the\ constant\ mean\ curvature\ and\ \psi_0\in\R}.
$$

{\bf case 2.} $\beta=1$ implies that the angle function which gives a CMC surface is an affine function
$\theta=2H x +d,\ {\rm with}\ H,d\in \R$.\\
For the general case of CMC hypersurfaces in $\mathbb{Q}_\epsilon^n\times\R$ see \cite[Thm. 3]{mn:To09}.
\end{remark}

\end{document}